\documentclass[11pt]{article}
\usepackage{lscape}
\usepackage{shadowtext}
\usepackage{tabularx,longtable}
\usepackage{tabularray}
\usepackage{comment}
\usepackage[most]{tcolorbox}
\usepackage{amsfonts}
\usepackage{multirow}
\usepackage{latexsym}
\usepackage{amssymb}
\usepackage{amsmath,blkarray,booktabs,bigstrut}
\usepackage{amsthm}
\usepackage{array}
\usepackage{setspace}
\usepackage{enumerate}
\usepackage{url}
\usepackage{bm}
\usepackage{tikz,color,lipsum,lmodern}
\usepackage{colortbl}
\usepackage{pgfplots}
\usepackage{graphicx}
\usepackage[export]{adjustbox}
\usepackage{caption}
\usepackage{subcaption}
\usepackage{tfrupee}
\usepackage{float}
\usepackage{fancyhdr}
\usepackage{centernot}
\usepackage{tkz-euclide}
\usepackage{tabto}
\usepackage{accents}
\usepackage{epsdice}

\usetikzlibrary{graphs, shapes.geometric, calc}
\usetikzlibrary{patterns, snakes, arrows.meta}
\usetikzlibrary{graphs,decorations.pathmorphing,decorations.markings}

\tikzstyle{startstop} = [rectangle, rounded corners, minimum width=3cm, minimum height=1cm,text centered, draw=black, fill=red!30]
\tikzstyle{decision} = [diamond, minimum width=3cm, minimum height=1cm, text centered, align=center, inner sep=0pt, draw=black, fill=green!30]
\tikzstyle{decision2} = [diamond, minimum width=3cm, minimum height=1cm, text centered, align=center, inner sep=0pt, draw=black, fill=green!30, aspect=2]
\tikzstyle{arrow} = [thick,->,>=stealth]
\tikzstyle{io} = [trapezium, trapezium left angle=70, trapezium right angle=110, minimum width=3cm, text width=3cm, minimum height=1cm, text centered, draw=black, fill=blue!30]
\tikzstyle{process} = [rectangle, minimum width=3cm, minimum height=1cm, text centered, draw=black, fill=orange!30]

\pgfmathsetmacro{\cellsize}{1cm}
\pgfmathsetmacro{\ladderstep}{\cellsize/6}
\pgfmathsetmacro{\ladderwidth}{\cellsize/6}
\tikzset{
	cell/.style={
		draw=green!50!black,top color=lime!10,bottom color=lime,
		line width=1pt,
		rounded corners=\cellsize/5
	},
	ladder/.style={
		draw=orange,line cap=round,line join=round,line width=1.5pt,
	},
	chute/.style={
		top color=blue!50,
		bottom color=blue!50,
		middle color=blue!5,
		draw=blue!75!black,line cap=round,
		line join=round,line width=0.5pt,
	},
	board auto/.style={
		x=\cellsize pt,
		y=\cellsize pt,
		decorate,decoration={
			markings,
			mark=between positions 0 and 1 step \cellsize pt with
			{\path[cell] (-.5*\cellsize pt,-.5*\cellsize pt)
				rectangle (.5*\cellsize pt,.5*\cellsize pt);}
		},
	},
	ladder draw/.style={
		decorate,decoration={
			markings,
			mark=between positions {1/#1/2} and {-1/#1/2} step {1/#1} with {
				\path[ladder]
				(0,-\ladderwidth pt) -- (0,\ladderwidth pt)
				(\ladderstep pt,-\ladderwidth pt) -- (-\ladderstep pt,-\ladderwidth pt)
				(\ladderstep pt,\ladderwidth pt) -- (-\ladderstep pt,\ladderwidth pt);
			},
		},
	},
	ladder auto/.style={
		x=\cellsize pt,
		y=\cellsize pt,
		to path={
			let
			\p1=($(\tikztostart) - (\tikztotarget)$),
			\n1={veclen(\x1,\y1)}
			in
			\pgfextra{
				\pgfmathsetmacro{\bars}{int(\n1/\ladderstep/2)+1}
				\pgfinterruptpath
				\draw[ladder draw=\bars] (\tikztostart) -- (\tikztotarget);
				\endpgfinterruptpath
			}
		},
	},
	ladder auto bis/.style={
		x=\cellsize pt,
		y=\cellsize pt,
		to path={
			let
			\p1=(\tikztostart),
			\p2=(\tikztotarget),
			\p3=($(\tikztostart) - (\tikztotarget)$),
			\p4=([xshift=-\ladderwidth pt]\p1),
			\p5=([xshift=-\ladderwidth pt]\p2),
			\p6=([xshift=\ladderwidth pt]\p1),
			\p7=([xshift=\ladderwidth pt]\p2),
			\n1={veclen(\x3,\y3)}
			in
			\pgfextra{
				\pgfmathtruncatemacro{\bars}{int(\n1/\ladderstep/2)+1}
				\pgfinterruptpath
				\path[ladder] (\p4) -- (\p5);
				\path[ladder] (\p6) -- (\p7)
				\foreach \bar in {1,...,\bars}{
					\pgfextra{\pgfmathsetmacro{\pos}{1/(\bars+1)*\bar}}
					coordinate[pos=\pos] (p\bar)
				};
				\foreach \bar in {1,...,\bars}{
					\path[ladder] (p\bar) -- ++(-2*\ladderwidth pt,0);
				}
				\endpgfinterruptpath
			}
		},
	},
	chute auto/.style={
		to path={
			let
			\p1=([xshift=\ladderwidth]\tikztostart),
			\p2=([xshift=-\ladderwidth]\tikztostart),
			\p3=([xshift=\ladderwidth]\tikztotarget),
			\p4=([xshift=-\ladderwidth]\tikztotarget),
			\p5=($(\p1)!.5!(\p3)$),
			\p6=($(\p2)!.5!(\p4)$)
			in
			\pgfextra{
				\pgfinterruptpath
				\path[thick,draw=blue]
				(\p1) sin (\p5) cos (\p3)
				(\p4) sin (\p6) cos (\p2);
				\fill[fill=cyan!20]
				(\p1) sin (\p5) cos (\p3) --
				(\p4) sin (\p6) cos (\p2) -- cycle;
				\endpgfinterruptpath
			}
		},
	},
	chute auto bis/.style={
		x=\cellsize pt,
		y=\cellsize pt,
		to path={
			let
			\p1=([xshift=\ladderwidth pt]\tikztostart),
			\p2=([xshift=-\ladderwidth pt]\tikztostart),
			\p3=([xshift=\ladderwidth pt]\tikztotarget),
			\p4=([xshift=-\ladderwidth pt]\tikztotarget),
			\p5=($(\p1)!.5!(\p3)$),
			\p6=($(\p2)!.5!(\p4)$)
			in
			\pgfextra{
				\pgfinterruptpath
				\path[chute]
				(\p1) sin (\p5) cos (\p3) --
				(\p4) sin (\p6) cos (\p2) -- cycle;
				\endpgfinterruptpath
			}
		},
	},
}
\tikzgraphsset{
	implies/.style={edges={double, double equal sign distance, -implies}}
}

\usepackage[superscript,nospace]{cite}
\usepackage{url}
\pgfplotsset{compat = newest}
\definecolor{grey}{rgb}{.5,.6,.5}

\newcommand{\p}{\mathrm{P}}

\newcommand{\bbN}{\mathbb{N}}
\newcommand{\cc}{^\mathsf{c}}
\DeclareMathOperator*{\inm}{\bar{\in}}
\DeclareMathOperator*{\nim}{\bar{\ni}}

\definecolor{grey}{rgb}{.7,.7,.7}
\definecolor{Blue}{rgb}{0.043,0.043,0.796}
\definecolor{dgrey}{rgb}{.4,.4,.4}
\definecolor{ggreen}{rgb}{0,.8,0}
\definecolor{cred}{rgb}{0.89,0.0,0.13}
\definecolor{lav}{rgb}{.4,.2,.6}	
\newtheorem{defn}{Definition}[section]
\newtheorem{result}{Result}[section]

\newtheorem{note}{Note}[section]
\theoremstyle{remark}

\pgfplotsset{
	standard/.style={set layers=tick labels on top},
	layers/tick labels on top/.define layer set=
	{axis background,
		axis grid,
		axis lines,
		main,
		axis tick labels,
		axis ticks,
		axis descriptions,
		axis foreground}
	{/pgfplots/layers/standard}
}

\begin{document}
	\title{\textbf{On the Game of Moksha-Patam}}
	\author{\textbf{Aninda Kumar Nanda$^1$ and Amit Kumar Misra$^2$}\\\textit{\small $^1$Indian Statistical Institute, Delhi Centre, New Delhi, India}\\\textit{\small $^2$Department of Statistics, Babasaheb Bhimrao Ambedkar University, Lucknow, India}}
	\date{\today}
	\maketitle 
	
	\begin{abstract}
		The game of \emph{Moksha-Patam}, often known as `Chutes and Ladders', is a widely played indoor game worldwide. While studies have been conducted regarding the nature of an individual board, the possibilities that open up when we change the positions of the chutes and the ladders on a board are mostly unventured. In this article, we classify and study the different possible types of \emph{Moksha-Patam} Boards, introduce a naming convention for them and establish a count for the boards that can be constructed adhering to some basic assumptions of the game.
	\end{abstract}

	\section{Introduction and Preliminaries}
	The game of `\emph{Moksha-Patam}' from India\cite{bierend} has been of interest to the subject of Markov Chains for a long time. Often alternatively known as `Chutes and Ladders' or `Snakes and Ladders', this game, as played on a particular design of the board with 10 chutes and 10 ladders, was studied by Daykin et. al.\cite{daykin} in 1967. In their study, Daykin et. al. explain in detail how the game can be studied as a Markov Chain. In doing so, they clarify the construction of the associated one-step probability transition matrix and how to obtain the probability distribution of the number of steps required to win the game. These techniques over the years have come to be the most commonly employed ones in the study of the game of `\emph{Moksha-Patam}'.\\
	
	In 1993, Altheon et. al.\cite{Althoen} introduced a particular `standard' board with 10 chutes and 9 ladders, and studied it using computer simulations besides the use of similar techniques as employed by Daykin et. al. However, unlike Daykin et. al. who studied the game as played from starting cell $1$, Altheon et. al. made the choice of starting the game from outside the board on some figurative cell $0$. Additionally, Altheon et. al. studied the impact of adding a single snake or ladder to the board being studied. In line with this thought, Diaconis and Durrett\cite{DD2000}, in 2000, studied the general effect of altering a single row of the one-step transition probability matrix of any Markov Chain. And while Altheon et. al. had once mentioned in the introduction of their work that the game played on their `standard' board is an absorbing Markov Chain, the one-step probability matrix of the board was closely studied by Cheteyan et. al.\cite{cheteyan} in 2011 and the claimed result was established mathematically. Besides studying this `standard' board using similar techniques, Hochman\cite{hochman} in 2014 concluded his REU work on the topic with certain generalised problems like the use of a general die rather than the regular six-faced fair die, and the effects of starting the game from an arbitrary cell on the board. While studying this `standard' board in 2019, Tae\cite{tae} made use of the idea of eigenvectors to obtain the stationary probability distributions of the Markov Chain. This led to yet another proof of the fact that the game played on this `standard' board is an absorbing Markov Chain.\\
	
	In 2021, Tun\cite{tun} used the original techniques established by Daykin et. al. to study the game as played on yet another board with 10 chutes and 10 ladders. The game played on this board, as Tun comments, is again an absorbing Markov Chain.\\
	
	While it may seem that the game, as played on any board, is an absorbing Markov Chain, a key takeaway from this study is that the game played on some particular boards may not be absorbing, i.e., the game may never end. In fact, some games will never end. We study the game of `\emph{Moksha-Patam}' by closely observing the one-step transition probability matrix of the Markov Chain, by establishing the existence of closed sets, stationary probability distributions and related ideas, and ultimately establish an algorithm to determine what boards lead to games that are absorbing Markov Chains and what boards lead to other different types of games. We also establish a count of the number of possible boards the game can be played on.\\

	We conduct the entirely of this study by adhering to the rules established by Daykin et. al. However, instead of studying one of the possible boards, here we shall study all possible boards on which the game of `\emph{Moksha-Patam}' can be played. In order to do so, we must establish the following few reasonable conventions:

	\begin{enumerate}
		\item The ladders and the chutes are essentially the same, in the sense that they are both one-sided wormholes between cells, their only difference being in the direction of the passage. Thus, it can be said, in a upward/downward directional sense, that a chute is an inverted ladder and vice-versa. Hence, we will, as required, refer to ladders and chutes as \emph{components} of the board, to avoid any notion of specific directionality. As and when the sense of directionality is required, we shall refer to the starting and ending cells of the components, keeping in mind that a ladder starts at the bottom and end at the top while a chute starts at the top and ends at the bottom.
		\item No two components on the board can start from the same cell [see Figure \ref{multiple_a}], for it would create confusion regarding which path to follow.
		\begin{figure}[H]
			\centering\subfloat[]{
				\begin{tikzpicture}
					\newcounter{c}
					\foreach \x in {0,1,2,3}
					\foreach \y in {7,8,9}
					\shade[inner color=white, outer color=grey] (\x,\y) rectangle (\x+1,\y+1);
					\draw[chute auto] (2.5,9.5) to (1.5,7.5);
					\draw[chute auto] (3.5,8.5) to (2.5,9.5);
					
					\setcounter{c}{81}
					\foreach \x in {1,2,3,4}
					\node at (\x-0.5,8+0.5) {\number\value{c}\addtocounter{c}{1}};
					\foreach \y in {7,9}
					\setcounter{c}{\y7}
					\foreach \x in {4,3,2,1}
					\node at (\x-0.5,\y+0.5) {\number\value{c}\addtocounter{c}{1}};
			\end{tikzpicture}}\hspace{0.35cm}
			\subfloat[]{
				\begin{tikzpicture}
					\foreach \x in {0,1,2,3}
					\foreach \y in {7,8,9}
					\shade[inner color=white, outer color=grey] (\x,\y) rectangle (\x+1,\y+1);
					\draw[chute auto] (3.5,8.5) to (1.5,7.5);
					\draw[ladder auto] (3.5,8.5) to (2.5,9.5);
					\setcounter{c}{81}
					\foreach \x in {1,2,3,4}
					\node at (\x-0.5,8+0.5) {\number\value{c}\addtocounter{c}{1}};
					\foreach \y in {7,9}
					\setcounter{c}{\y7}
					\foreach \x in {4,3,2,1}
					\node at (\x-0.5,\y+0.5) {\number\value{c}\addtocounter{c}{1}};
			\end{tikzpicture}}\hspace{0.35cm}
			\subfloat[]{
				\begin{tikzpicture}
					\foreach \x in {0,1,2,3}
					\foreach \y in {7,8,9}
					\shade[inner color=white, outer color=grey] (\x,\y) rectangle (\x+1,\y+1);
					\draw[ladder auto] (2.5,9.5) to (1.5,7.5);
					\draw[ladder auto] (3.5,8.5) to (1.5,7.5);
					\setcounter{c}{81}
					\foreach \x in {1,2,3,4}
					\node at (\x-0.5,8+0.5) {\number\value{c}\addtocounter{c}{1}};
					\foreach \y in {7,9}
					\setcounter{c}{\y7}
					\foreach \x in {4,3,2,1}
					\node at (\x-0.5,\y+0.5) {\number\value{c}\addtocounter{c}{1}};
			\end{tikzpicture}}
			\caption{Components starting at the  same cell is not allowed.}\label{multiple_a}
		\end{figure}
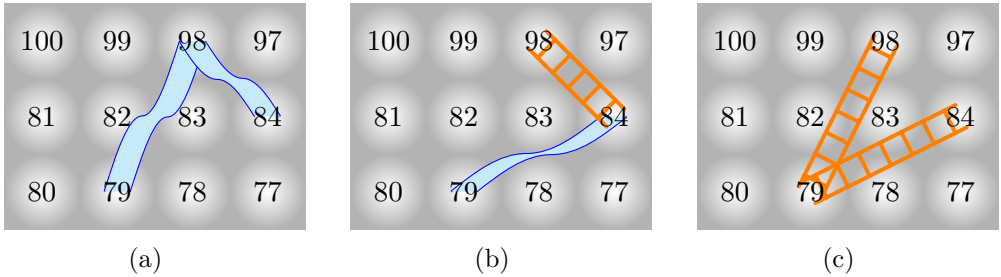
		If one component starts from the cell where another component ends, then the same overall effect is obtained by replacing the latter by a component that has the same starting cell but ends at the ending of the former [see Figure \ref{multiple_b}]. Now, two components can share an ending cell, and nothing in the rules of the game suggest against such a scenario. Throughout most of our study we consider all possible boards, including those where components share their ending cells.
		\item The player starts the game from cell $1$ and reaching the cell $100$ marks the completion of a game. Hence, having a ladder start at cell $1$, doesn't make much sense since it will imply availing the ladder before the game even begins. Likewise, having a chute start at cell $100$, doesn't make much sense since it will imply availing the chute after the game has ended. While there is no harm in ending a chute on cell $1$, it just means all progress of the game is completely lost. Similarly, having a ladder end on cell $100$ means there is some cell other than $100$ which also acts as the winning cell, i.e., there is a cell other than $100$ reaching which we will win. To avoid such duality of the ultimate goal, and to avoid essentially restarting the game while playing, we shall adhere to the norm that cell $1$ and cell $100$ cannot house any component. 
		
		\begin{figure}[H]
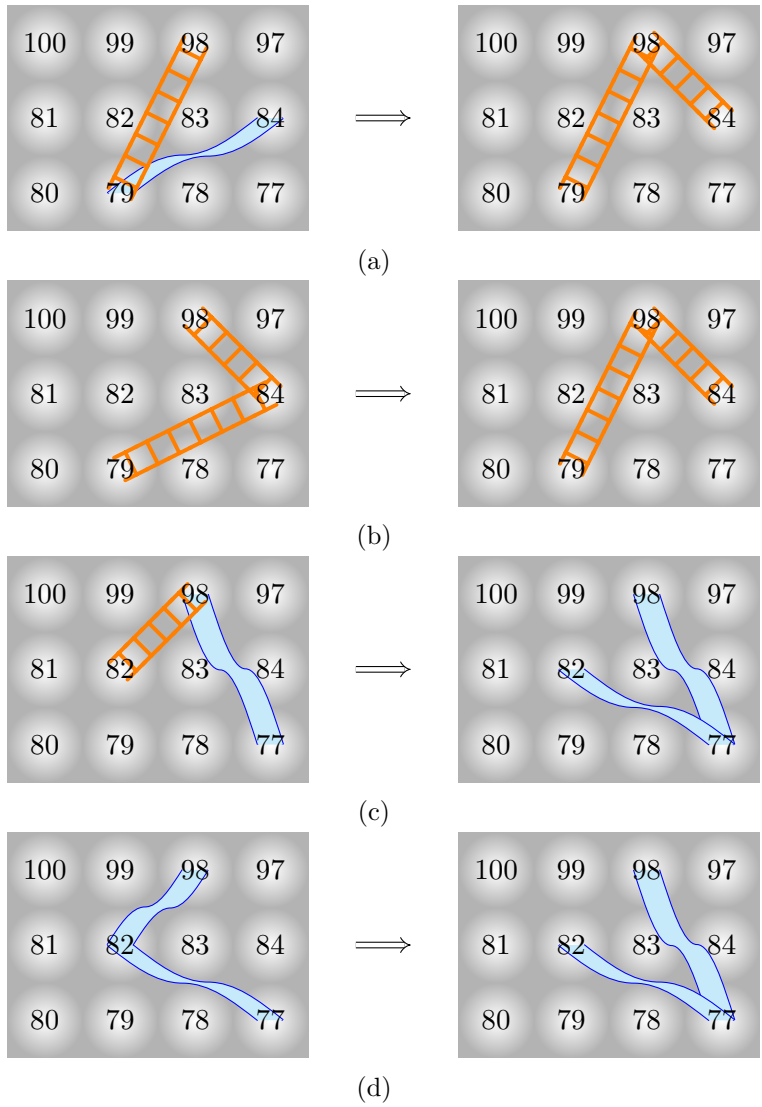

			\centering\subfloat[]{
				\tikz[>=implies]{
					\foreach \x in {0,1,2,3}
					\foreach \y in {7,8,9}
					\shade[inner color=white, outer color=grey] (\x,\y) rectangle (\x+1,\y+1);
					\draw[chute auto] (3.5,8.5) to (1.5,7.5);
					\draw[ladder auto] (1.5,7.5) to (2.5,9.5);
					\setcounter{c}{81}
					\foreach \x in {1,2,3,4}
					\node at (\x-0.5,8+0.5) {\number\value{c}\addtocounter{c}{1}};
					\foreach \y in {7,9}
					\setcounter{c}{\y7}
					\foreach \x in {4,3,2,1}
					\node at (\x-0.5,\y+0.5) {\number\value{c}\addtocounter{c}{1}};
					
					\node (istart) at (4.5,8.5) {};
					\node (istop) at (5.5,8.5) {};
					\graph[implies] {(istart) -> (istop)};
					
					\foreach \x in {6,7,8,9}
					\foreach \y in {7,8,9}
					\shade[inner color=white, outer color=grey] (\x,\y) rectangle (\x+1,\y+1);
					\draw[ladder auto] (9.5,8.5) to (8.5,9.5);
					\draw[ladder auto] (7.5,7.5) to (8.5,9.5);
					\setcounter{c}{81}
					\foreach \x in {7,8,9,10}
					\node at (\x-0.5,8+0.5) {\number\value{c}\addtocounter{c}{1}};
					\foreach \y in {7,9}
					\setcounter{c}{\y7}
					\foreach \x in {10,9,8,7}
					\node at (\x-0.5,\y+0.5) {\number\value{c}\addtocounter{c}{1}};
			}}\\\subfloat[]{
				\tikz[>=implies]{
					\foreach \x in {0,1,2,3}
					\foreach \y in {7,8,9}
					\shade[inner color=white, outer color=grey] (\x,\y) rectangle (\x+1,\y+1);
					\draw[ladder auto] (3.5,8.5) to (1.5,7.5);
					\draw[ladder auto] (3.5,8.5) to (2.5,9.5);
					\setcounter{c}{81}
					\foreach \x in {1,2,3,4}
					\node at (\x-0.5,8+0.5) {\number\value{c}\addtocounter{c}{1}};
					\foreach \y in {7,9}
					\setcounter{c}{\y7}
					\foreach \x in {4,3,2,1}
					\node at (\x-0.5,\y+0.5) {\number\value{c}\addtocounter{c}{1}};
					
					\node (istart) at (4.5,8.5) {};
					\node (istop) at (5.5,8.5) {};
					\graph[implies] {(istart) -> (istop)};
					
					\foreach \x in {6,7,8,9}
					\foreach \y in {7,8,9}
					\shade[inner color=white, outer color=grey] (\x,\y) rectangle (\x+1,\y+1);
					\draw[ladder auto] (9.5,8.5) to (8.5,9.5);
					\draw[ladder auto] (7.5,7.5) to (8.5,9.5);
					\setcounter{c}{81}
					\foreach \x in {7,8,9,10}
					\node at (\x-0.5,8+0.5) {\number\value{c}\addtocounter{c}{1}};
					\foreach \y in {7,9}
					\setcounter{c}{\y7}
					\foreach \x in {10,9,8,7}
					\node at (\x-0.5,\y+0.5) {\number\value{c}\addtocounter{c}{1}};
			}}\\
			
			\subfloat[]{
				\tikz[>=implies]{
					\foreach \x in {0,1,2,3}
					\foreach \y in {7,8,9}
					\shade[inner color=white, outer color=grey] (\x,\y) rectangle (\x+1,\y+1);
					\draw[chute auto] (3.5,7.5) to (2.5,9.5);
					\draw[ladder auto] (1.5,8.5) to (2.5,9.5);
					\setcounter{c}{81}
					\foreach \x in {1,2,3,4}
					\node at (\x-0.5,8+0.5) {\number\value{c}\addtocounter{c}{1}};
					\foreach \y in {7,9}
					\setcounter{c}{\y7}
					\foreach \x in {4,3,2,1}
					\node at (\x-0.5,\y+0.5) {\number\value{c}\addtocounter{c}{1}};
					
					\node (istart) at (4.5,8.5) {};
					\node (istop) at (5.5,8.5) {};
					\graph[implies] {(istart) -> (istop)};
					
					\foreach \x in {6,7,8,9}
					\foreach \y in {7,8,9}
					\shade[inner color=white, outer color=grey] (\x,\y) rectangle (\x+1,\y+1);
					\draw[chute auto] (9.5,7.5) to (8.5,9.5);
					\draw[chute auto] (7.5,8.5) to (9.5,7.5);
					\setcounter{c}{81}
					\foreach \x in {7,8,9,10}
					\node at (\x-0.5,8+0.5) {\number\value{c}\addtocounter{c}{1}};
					\foreach \y in {7,9}
					\setcounter{c}{\y7}
					\foreach \x in {10,9,8,7}
					\node at (\x-0.5,\y+0.5) {\number\value{c}\addtocounter{c}{1}};
			}}\\\subfloat[]{
				\tikz[>=implies]{
					\foreach \x in {0,1,2,3}
					\foreach \y in {7,8,9}
					\shade[inner color=white, outer color=grey] (\x,\y) rectangle (\x+1,\y+1);
					\draw[chute auto] (3.5,7.5) to (1.5,8.5);
					\draw[chute auto] (1.5,8.5) to (2.5,9.5);
					\setcounter{c}{81}
					\foreach \x in {1,2,3,4}
					\node at (\x-0.5,8+0.5) {\number\value{c}\addtocounter{c}{1}};
					\foreach \y in {7,9}
					\setcounter{c}{\y7}
					\foreach \x in {4,3,2,1}
					\node at (\x-0.5,\y+0.5) {\number\value{c}\addtocounter{c}{1}};
					
					\node (istart) at (4.5,8.5) {};
					\node (istop) at (5.5,8.5) {};
					\graph[implies] {(istart) -> (istop)};
					
					\foreach \x in {6,7,8,9}
					\foreach \y in {7,8,9}
					\shade[inner color=white, outer color=grey] (\x,\y) rectangle (\x+1,\y+1);
					\draw[chute auto] (9.5,7.5) to (8.5,9.5);
					\draw[chute auto] (7.5,8.5) to (9.5,7.5);
					\setcounter{c}{81}
					\foreach \x in {7,8,9,10}
					\node at (\x-0.5,8+0.5) {\number\value{c}\addtocounter{c}{1}};
					\foreach \y in {7,9}
					\setcounter{c}{\y7}
					\foreach \x in {10,9,8,7}
					\node at (\x-0.5,\y+0.5) {\number\value{c}\addtocounter{c}{1}};
			}}
			\caption{Components sharing a starting cell and an ending cell is equivalent to components sharing their ending cells.}\label{multiple_b}
		\end{figure}
		
	\end{enumerate}


	Besides the conventions listed above, we shall make use of some specific terminologies as defined below:
	
	\begin{defn}[Entrance]
		A cell on a \emph{Moksha-Patam} board is said to be an \emph{entrance} if a component starts on that cell, i.e., the cell houses either the foot of a ladder or the top of a chute.
	\end{defn}
	
	\begin{defn}[Non-entrance Cell]
		A cell which is not an entrance is said to be a \emph{non-entrance cell}.
	\end{defn}
	
	\begin{defn}[Exit]
		A cell on a \emph{Moksha-Patam} board is said to be an \emph{exit} if a component ends on that cell, i.e., the cell houses either the peak of a ladder or the bottom of a chute.
	\end{defn}
	
	\begin{defn}[Non-exit Cell]
		A cell which is not an exit is said to be a \emph{non-exit cell}.
	\end{defn}
	
	\begin{defn}[Non-component Cell]
		A cell on a \emph{Moksha-Patam} board is said to be a \emph{non-component cell} if it is neither an entrance nor an exit.
	\end{defn}
	
	
	\begin{defn}[Chute-Barrier]\label{chutebarr}
		A collection of six or more chutes is called a \emph{chute-barrier} if they have consecutive entrances.
	\end{defn}

	Now, to distinguish between all possible distinct boards from each other, we introduce a naming convention.
	
	\section{A Naming Convention}

	The naming convention we shall be using henceforth is as follows:
	\begin{quote}
		A board shall be called an ``$N(X)$ Board''  if it contains $N$ components (either ladders or chutes) on the board, with $X$ being a $2\times N$ matrix, the $(1,i)^\text{th}$ and $(2,i)^\text{th}$ elements of which specify the starting and ending positions of the $i^\text{th}$ component, respectively.
	\end{quote}
	
	To facilitate the uniqueness of the name of a particular board, we shall always list the starting positions of the components of the board, in $X$, in ascending order. That is the name ``$2\left(\begin{bmatrix}23&10\\5&60\end{bmatrix}\right)$ Board'' is not an acceptable name. Rather, the correct name, according to our convention, is ``$2\left(\begin{bmatrix}10&23\\60&5\end{bmatrix}\right)$ Board''.
	
	It is worth noting that an $N\bigg(((k_{ij}))\bigg)$ Board has a chute (or, ladder) if for some $j\in\bbN_N$, $k_{1j}>(\text{or,}<)\,k_{2j}$.\\
	
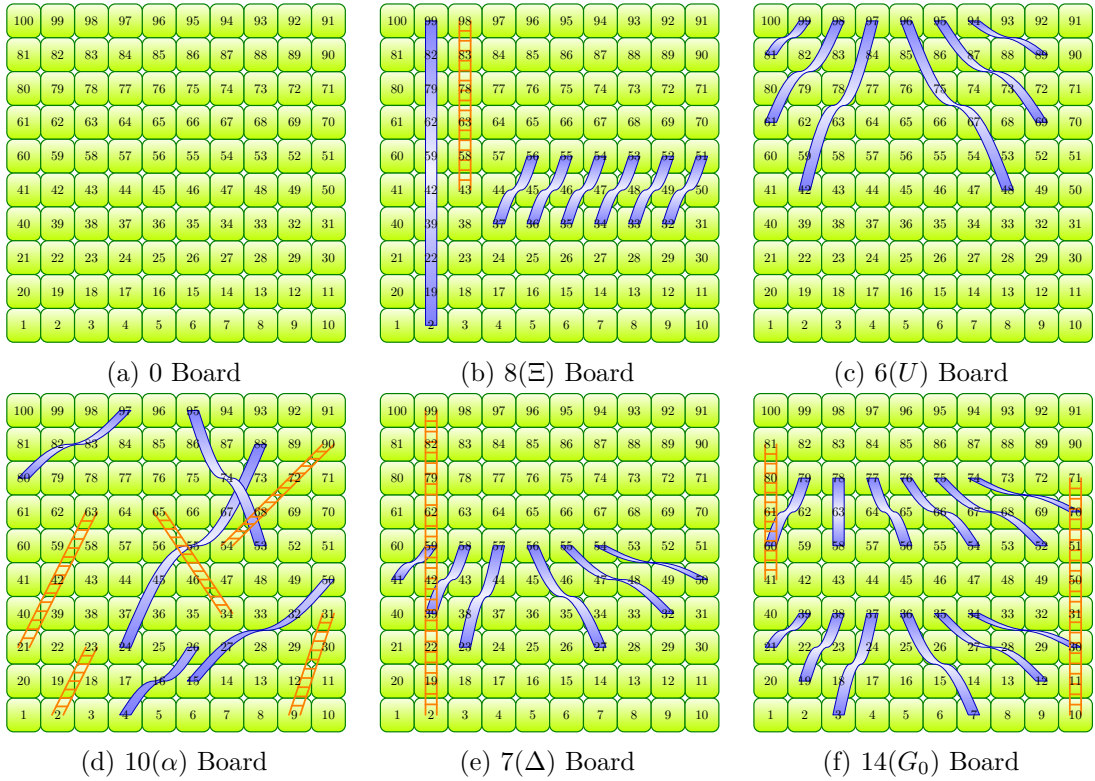
\begin{figure}[H]
		\centering
		\subfloat[$0$ Board\label{empty}]{\scalebox{0.45}{
		\begin{tikzpicture}	
			\draw[board auto] (1,1) -- (10,1) -- (10,10) -- (1,10) -- (1,2) -- (9,2) -- (9,9) -- (2,9) -- (2,3) -- (8,3) -- (8,8) -- (3,8) -- (3,4) -- (7,4) -- (7,7) -- (4,7) -- (4,5) -- (6,5) -- (6,6) -- (5,6);	
			
			\foreach \y in {0,2,4,6,8}
			\setcounter{c}{\y1}
			\foreach \x in {1,2,3,4,5,6,7,8,9,10}
			\node at (\x,\y+1) {\number\value{c}\addtocounter{c}{1}};
			\foreach \y in {1,3,5,7,9}
			\setcounter{c}{\y1}
			\foreach \x in {10,9,8,7,6,5,4,3,2,1}
			\node at (\x,\y+1) {\number\value{c}\addtocounter{c}{1}};
		\end{tikzpicture}}}\hspace{0.2cm}
		\subfloat[$8(\Xi)$ Board\label{8Xi}]{\scalebox{0.45}{
		\begin{tikzpicture}	
			\draw[board auto] (1,1) -- (10,1) -- (10,10) -- (1,10) -- (1,2) -- (9,2) -- (9,9) -- (2,9) -- (2,3) -- (8,3) -- (8,8) -- (3,8) -- (3,4) -- (7,4) -- (7,7) -- (4,7) -- (4,5) -- (6,5) -- (6,6) -- (5,6);	
			\draw[chute auto bis] (10,6) to (9,4) (9,6) to (8,4) (8,6) to (7,4) (7,6) to (6,4) (6,6) to (5,4) (5,6) to (4,4) (2,10) to (2,1);
			\draw[ladder auto bis] (3,5) to (3,10);
			
			\foreach \y in {0,2,4,6,8}
			\setcounter{c}{\y1}
			\foreach \x in {1,2,3,4,5,6,7,8,9,10}
			\node at (\x,\y+1) {\number\value{c}\addtocounter{c}{1}};
			\foreach \y in {1,3,5,7,9}
			\setcounter{c}{\y1}
			\foreach \x in {10,9,8,7,6,5,4,3,2,1}
			\node at (\x,\y+1) {\number\value{c}\addtocounter{c}{1}};
		\end{tikzpicture}}
		}\hspace{0.2cm}
		\subfloat[$6(U)$ Board\label{6U}]{\scalebox{0.45}{
		\begin{tikzpicture}	
			\draw[board auto] (1,1) -- (10,1) -- (10,10) -- (1,10) -- (1,2) -- (9,2) -- (9,9) -- (2,9) -- (2,3) -- (8,3) -- (8,8) -- (3,8) -- (3,4) -- (7,4) -- (7,7) -- (4,7) -- (4,5) -- (6,5) -- (6,6) -- (5,6);	
			\draw[chute auto bis] (2,10) to (1,9) (3,10) to (1,7) (4,10) to (2,5) (5,10) to (8,5) (6,10) to (9,7) (7,10) to (9,9);
			
			\foreach \y in {0,2,4,6,8}
			\setcounter{c}{\y1}
			\foreach \x in {1,2,3,4,5,6,7,8,9,10}
			\node at (\x,\y+1) {\number\value{c}\addtocounter{c}{1}};
			\foreach \y in {1,3,5,7,9}
			\setcounter{c}{\y1}
			\foreach \x in {10,9,8,7,6,5,4,3,2,1}
			\node at (\x,\y+1) {\number\value{c}\addtocounter{c}{1}};
		\end{tikzpicture}}}\\
		\subfloat[$10(\alpha)$ Board\label{10alpha}]{\scalebox{0.45}{
				\begin{tikzpicture}	
					\draw[board auto] (1,1) -- (10,1) -- (10,10) -- (1,10) -- (1,2) -- (9,2) -- (9,9) -- (2,9) -- (2,3) -- (8,3) -- (8,8) -- (3,8) -- (3,4) -- (7,4) -- (7,7) -- (4,7) -- (4,5) -- (6,5) -- (6,6) -- (5,6);	
					\draw[chute auto bis] (8,9) to (4,3) (4,10) to (1,8) (6,10) to (8,6) (10,5) to (6,2) (6,3) to (4,1);
					\draw[ladder auto bis] (2,1) to (3,3) (1,3) to (3,7) (7,6) to (10,9) (7,4) to (5,7) (9,1) to (10,4);
					
					\foreach \y in {0,2,4,6,8}
					\setcounter{c}{\y1}
					\foreach \x in {1,2,3,4,5,6,7,8,9,10}
					\node at (\x,\y+1) {\number\value{c}\addtocounter{c}{1}};
					\foreach \y in {1,3,5,7,9}
					\setcounter{c}{\y1}
					\foreach \x in {10,9,8,7,6,5,4,3,2,1}
					\node at (\x,\y+1) {\number\value{c}\addtocounter{c}{1}};
		\end{tikzpicture}}}\hspace{0.2cm}
		\subfloat[$7(\Delta)$ Board\label{7Delta}]{\scalebox{0.45}{
				\begin{tikzpicture}	
					\draw[board auto] (1,1) -- (10,1) -- (10,10) -- (1,10) -- (1,2) -- (9,2) -- (9,9) -- (2,9) -- (2,3) -- (8,3) -- (8,8) -- (3,8) -- (3,4) -- (7,4) -- (7,7) -- (4,7) -- (4,5) -- (6,5) -- (6,6) -- (5,6);	
					\draw[chute auto bis] (2,6) to (1,5) (3,6) to (2,4) (4,6) to (3,3) (5,6) to (7,3) (6,6) to (9,4) (7,6) to (10,5);
					\draw[ladder auto bis] (2,1) to (2,10);
					
					\foreach \y in {0,2,4,6,8}
					\setcounter{c}{\y1}
					\foreach \x in {1,2,3,4,5,6,7,8,9,10}
					\node at (\x,\y+1) {\number\value{c}\addtocounter{c}{1}};
					\foreach \y in {1,3,5,7,9}
					\setcounter{c}{\y1}
					\foreach \x in {10,9,8,7,6,5,4,3,2,1}
					\node at (\x,\y+1) {\number\value{c}\addtocounter{c}{1}};
			\end{tikzpicture}}
		}\hspace{0.2cm}
		\subfloat[$14(G_0)$ Board\label{14G_0}]{\scalebox{0.45}{
				\begin{tikzpicture}	
					\draw[board auto] (1,1) -- (10,1) -- (10,10) -- (1,10) -- (1,2) -- (9,2) -- (9,9) -- (2,9) -- (2,3) -- (8,3) -- (8,8) -- (3,8) -- (3,4) -- (7,4) -- (7,7) -- (4,7) -- (4,5) -- (6,5) -- (6,6) -- (5,6);	
					\draw[chute auto bis] (2,4) to (1,3) (3,4) to (2,2) (4,4) to (3,1) (5,4) to (7,1) (6,4) to (9,2) (7,4) to (10,3);
					\draw[chute auto bis] (2,8) to (1,6) (3,8) to (3,6) (4,8) to (5,6) (5,8) to (7,6) (6,8) to (9,6) (7,8) to (10,7);
					\draw[ladder auto bis] (10,1) to (10,8) (1,5) to (1,9);
					
					\foreach \y in {0,2,4,6,8}
					\setcounter{c}{\y1}
					\foreach \x in {1,2,3,4,5,6,7,8,9,10}
					\node at (\x,\y+1) {\number\value{c}\addtocounter{c}{1}};
					\foreach \y in {1,3,5,7,9}
					\setcounter{c}{\y1}
					\foreach \x in {10,9,8,7,6,5,4,3,2,1}
					\node at (\x,\y+1) {\number\value{c}\addtocounter{c}{1}};
		\end{tikzpicture}}}
		\caption{A Few \emph{Moksha-Patam} Boards.}
		\label{listboards}
	\end{figure}

	In Figure \ref{listboards}, we list a few boards along with their appropriate names as per our naming convention. Here, the matrices $\Xi$, $U$, $\alpha$, $\Delta$ and $G_0$ are given by $$\Xi=\begin{bmatrix}43&51&52&53&54&55&56&99\\98&32&33&34&35&36&37&2\end{bmatrix},$$ $$U=\begin{bmatrix}94&95&96&97&98&99\\89&69&48&42&61&81\end{bmatrix},$$ $$\alpha=\begin{bmatrix}2&9&21&26&34&50&54&88&95&97\\ 23&31&63&4&65&15&90&24&53&80\end{bmatrix},$$
	$$\Delta=\begin{bmatrix}2&54&55&56&57&58&59\\99&50&32&27&23&39&41\end{bmatrix}, \text{and}$$
	$$G_0=\begin{bmatrix}10&34&35&36&37&38&39&41&74&75&76&77&78&79\\71&30&12&7&3&19&21&81&70&52&54&56&58&60\end{bmatrix}.$$

\section{Classification of Boards}

	Note that the $0$ Board can always be completed, if the game is played for long enough, whereas the $6(U)$ Board can never be completed. The $7(\Delta)$ Board, however, can be completed if and only if the first die roll yields {\large\epsdice{1}}, that is the game can be won in a finite number of moves with a probability of $\frac16$.\\
	
	Hence, we classify the \emph{Moksha-Patam} boards into three categories:
	\begin{itemize}
		\item \textbf{Unwinnable Boards:} An example of such a board would be the $6(U)$ Board, shown in Figure \ref{6U}. These are boards that can never be completed, i.e., cell $100$ can never be reached.
		\item \textbf{Occasionally Winnable Boards:} An example of such a board would be the $7(\Delta)$ Board, shown in Figure \ref{7Delta}. These are boards that have a positive probability for both, never completing the game, and completing the game in a finite number of moves.
		\item \textbf{Ultimately Winnable Boards:} An example of such a board would be the $0$ Board, shown in Figure \ref{empty}. These are the \emph{Moksha-Patam} boards which will ultimately be completed if one plays long enough.
	\end{itemize}

	We can mathematically formalize these three categories by studying the game of \emph{Moksha-Patam} played on the $N(X)$ Board as a Markov Chain, $\{^{N(X)}M_n\}$, with state space $\bbN_{100}=\{n\in\bbN:n\leq 100\}$. We may drop the superscript specifying the board if the particular board the game is played on is clear from the context. In order to formalize the mathematical criteria of the classifications, we'll make use of the following results, a detailed discussion upon which can be found in Section 3.9 of Medhi[1982]\cite{medhi}.

\begin{result}\label{upzero}
	The one-step transition probability matrix $P$ of a discrete time-homogenous Markov Chain with a conveniently rearranged finite state space $S$, can be expressed as $$P_{|S|\times|S|}=\begin{bmatrix}P_1&\mathbb{O}_{|S_1|\times|S_2|}\\Q&P_2\end{bmatrix},$$ with $P_1$ being the $|S_1|\times|S_1|$ matrix of one-step transition probabilities between the states in $S_1\subset S$, $\mathbb{O}$ being the null matrix of order $|S_1|\times|S_2|$, $Q$ being the $|S_2|\times|S_1|$ matrix of one-step transition probabilities from the states in $S_2=S\backslash S_1$ to the states in $S_1$, and $P_2$ being the $|S_2|\times|S_2|$ matrix of one-step transition probabilities between the states in $S_2$ if, and only if, $\forall i\in S_1$ and $\forall j\in S_2$, $i\centernot\rightarrow j,$ that is, $S_1$ is a closed set.
\end{result}

\begin{result}\label{downzero}
	The one-step transition probability matrix $P$ of a discrete time-homogenous Markov Chain with a conveniently rearranged finite state space $S$, can be expressed as $$P_{|S|\times|S|}=\begin{bmatrix}P_1&Q\\\mathbb{O}_{|S_2|\times|S_1|}&P_2\end{bmatrix},$$ with $P_1$ being the $|S_1|\times|S_1|$ matrix of one-step transition probabilities between the states in $S_1\subset S$, $Q$ being the $|S_1|\times|S_2|$ matrix of one-step transition probabilities from the states in $S_1$ to the states in $S_2=S\backslash S_1$, $\mathbb{O}$ being the null matrix of order $|S_2|\times|S_1|$, and $P_2$ being the $|S_2|\times|S_2|$ matrix of one-step transition probabilities between the states in $S_2$ if, and only if, $\forall i\in S_1$ and $\forall j\in S_2$, $j\centernot\rightarrow i$, that is, $S_2$ is a closed set.
\end{result}

	Now, suppose that a board is unwinnable. That is $1\centernot\rightarrow100$. But since accessibility is a transitive relation, $\forall k\in\bbN_{100}$ such that $1\rightarrow k$, we must have $k\centernot\rightarrow100$. Likewise, $\forall k\in\bbN_{100}$ such that $k\rightarrow100$, we must have $1\centernot\rightarrow k$. Hence, we can partition $\bbN_{100}$ into two subsets $S_1$ and $S_{100}$ such that $1\in S_1$, $100\in S_{100}$, and $\forall i \in S_1$ and $\forall j\in S_{100}$, $i\centernot\rightarrow j$. Of course, there can be some state $\ell\in\bbN_{100}$ such that $1\centernot\rightarrow\ell\centernot\rightarrow100$, and we can include $\ell$ in either of $S_1$ and $S_{100}$.\\
	
	Thus, by Result \ref{upzero}, rearranging the states in $\bbN_{100}$ such that the states in $S_1$ are listed before the states in $S_{100}$, we must get a one-step transition probability matrix that can be written as $$P_{(r)}=\begin{bmatrix}P_1&\mathbb{O}\\Q&P_{100}\end{bmatrix},$$ where we use the subscript `$(r)$' to denote a rearranged one-step transition probability matrix. It is worth noting that such a rearrangement is done by pairs of consecutive elementary matrix operations. For example, to shift the $i^\text{th}$ state to the $j^\text{th}$ position, we first shift the $i^\text{th}$ row of the matrix to the $j^\text{th}$ row and then shift the $i^\text{th}$ column of the matrix to the $j^\text{th}$ column. \\
	
	Hence, we have the following necessary and sufficient criteria of an unwinnable board:
	\begin{result}[The Necessary and Sufficient Condition for Unwinnable Boards]\label{unwin}
		A board is unwinnable if, and only if, its associated one-step transition probability matrix $P$ can be rewritten as $$P_{(r)}=\begin{bmatrix}P_1&\mathbb{O}\\Q&P_{100}\end{bmatrix},$$ with square matrices $P_1\,(\nim p_{1,1})$ and $P_{100}\,(\nim p_{100,100})$, by rearranging the order in which the states are listed.
	\end{result}Here, we use the notation $p_{k,l}\inm P_*$ to indicate that the matrix $P_*$ contains the element $p_{k,l}=\p[M_1=l|M_0=k]$, where .\\

\begin{figure}[H]
		\centering
		\subfloat[$P$]{\scalebox{0.97}{
				\begin{tikzpicture}
					\node (fig) at (0,0) {\includegraphics[scale=0.3]{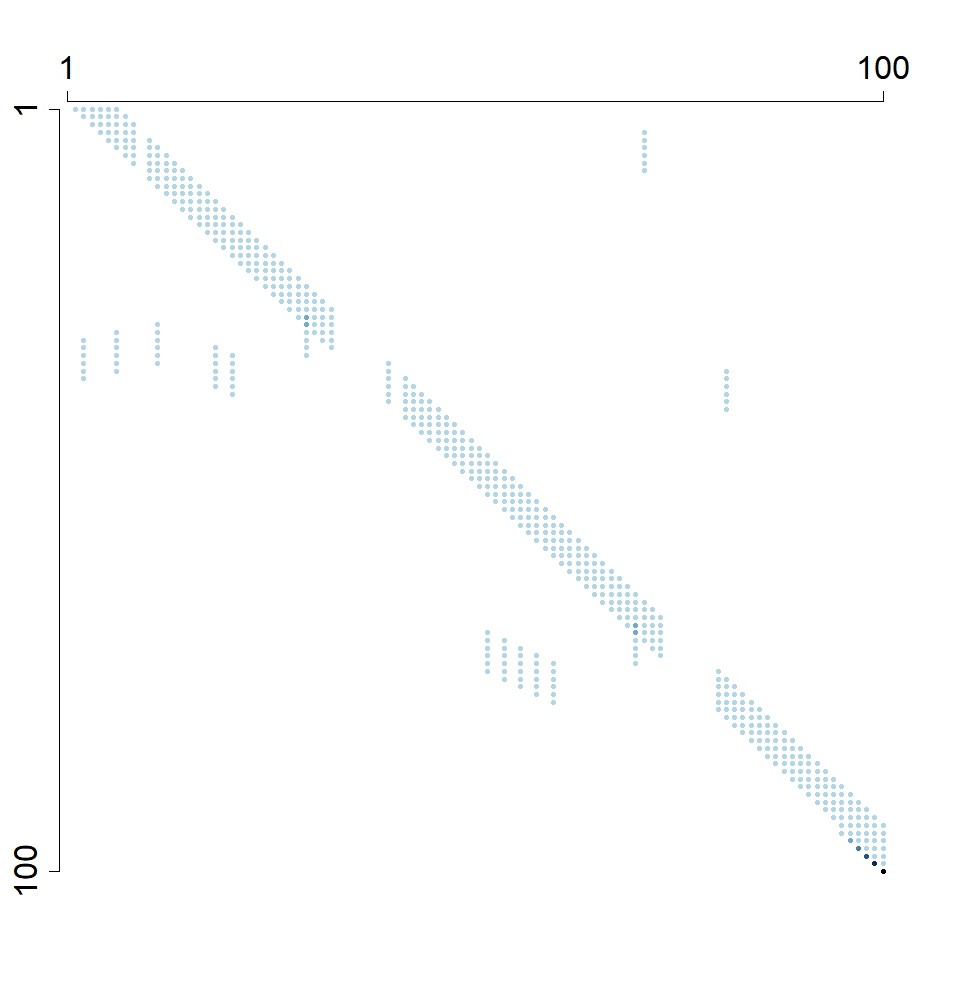}};
		\end{tikzpicture}}}
		\subfloat[$P_{(r)}$]{\scalebox{0.97}{
				\begin{tikzpicture}
					\node (fig) at (0,0) {\includegraphics[scale=0.3]{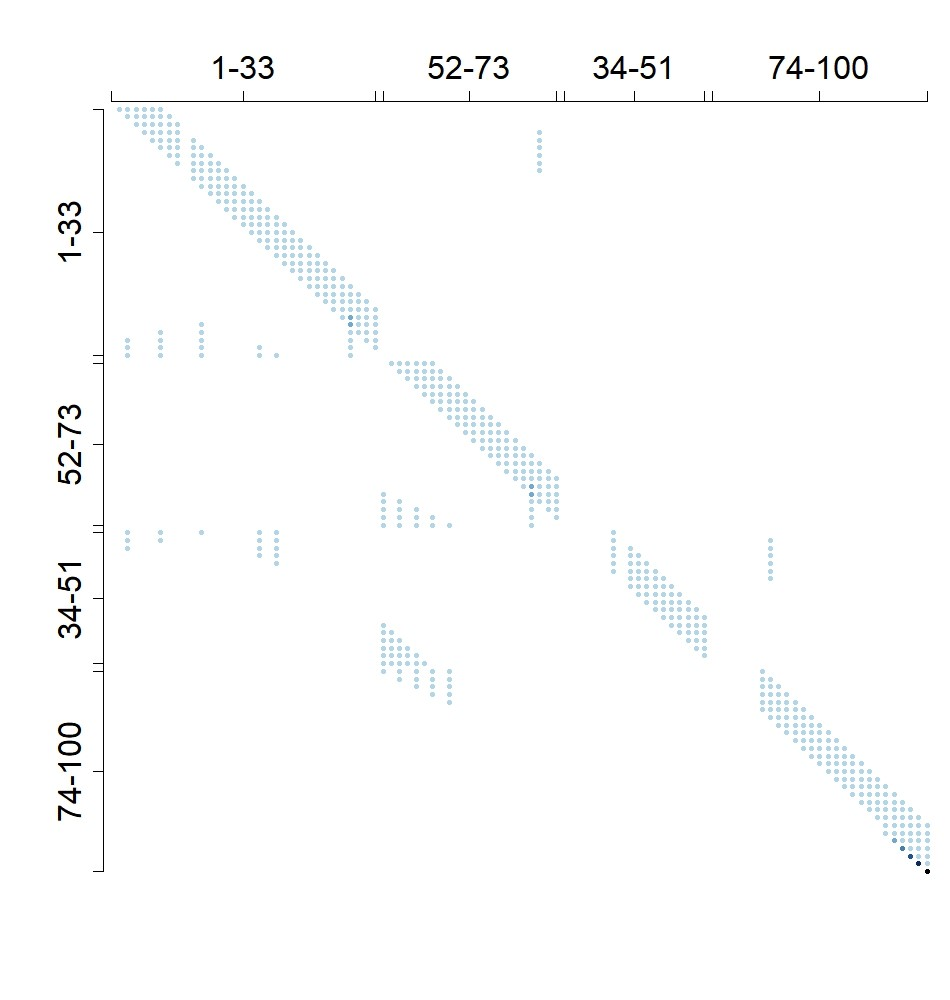}};
					\draw[cred, dashed] (-2.8,-0.24) -- (3.5,-0.24);
					\draw[cred, dashed] (0.693,-3) -- (0.693,3);
			\end{tikzpicture}}	
		}
		\caption{One-step Transition Probability Matrix of the $14(G_0)$ Board.}\label{14G0r}
	\end{figure}

	Figure \ref{14G0r} compares the arrangement of the typical one-step transition probability matrix of the $14(G_0)$ Board, with that of its rearranged one-step transition probability matrix. Here, each probability is indicated by an appropriate shade between white (indicating 0) and black (indicating 1). It is evident that the right-upper sub-matrix of $P_{(r)}$ is a null matrix which clearly exhibits, as stated in Result \ref{unwin}, that the $14(G_0)$ Board is unwinnable.
	
	\begin{result}[The Necessary and Sufficient Condition for Occasionally Winnable Boards]\label{occwin}
		A winnable (i.e., not unwinnable) board is occasionally winnable if, and only if,
		\begin{enumerate}[(i)]
			\item the associated one-step transition probability matrix $P$ can be rewritten as $$P_{(r)}=\begin{bmatrix}P_0&Q\\\mathbb{O}&P_C\end{bmatrix},$$ with some square matrices $P_0\,(\nim p_{1,100})$ and $P_C$, by rearranging the order in which the states are listed, and
			\item for every such rearrangement satisfying the above condition, $$P_{(r)}\neq\begin{bmatrix}P_0^*&\mathbb{O}\\Q^*&P_{C_1}\end{bmatrix}$$ for any square matrices $P_0^*\,(\nim p_{1,100})$ and $P_{C_1}\,(\rhd P_C)$.
		\end{enumerate}
	\end{result} Here, we use the notation $A\lhd B$ to indicate that $A$ is a sub-matrix of $B$.
	
	\begin{proof}
		\textit{Sufficient Part:}\\
		Let $N(X)$ Board be winnable (not unwinnable) and have the two properties listed above. Let $P_C$ contain the one-step transition probabilities amongst the states belonging to the non-empty set $C\subset \bbN_{100}$. Likewise, $P_0$ contains the one-step transition probabilities amongst the states belonging to $C\cc=\bbN_{100}\backslash C$. Hence, by Result \ref{downzero}, $C$ is a closed set. Also, since $p_{1,100}\inm P_0$, we have $1,100\in C\cc$. That is, $1,100\notin C$. \\
		
		However, since $P_{(r)}\neq\begin{bmatrix}P_0^*&\mathbb{O}\\Q^*&P_{C_1}\end{bmatrix}$ with $p_{1,100}\inm P_0^*$ and $P_C\lhd P_{C_1}$, no superset of $\{1,100\}$, which is also a subset of $C\cc$, is a closed set [see Result \ref{upzero}].\\
		
		Thus, $\{1,100\}$ is not a closed set, implying either $\exists\,j\in C$ such that $1\rightarrow j$, or $\exists\, k_1\in C\cc\backslash\{1,100\}$ such that $1\rightarrow k_1$. In the instance that the second case be correct, note that $\{1,k_1,100\}$ is not a closed set, implying either $\exists\, j\in C$ such that $1\rightarrow j$, or $\exists\, k_2\in C\cc\backslash\{1,k_1,100\}$ such that $1\rightarrow k_2$. By the Principle of Mathematical Induction, we can conclude that either  $\exists\, j\in C$ such that $1\rightarrow j$ or $\exists\,k_\eta\in C\cc\backslash\{1,k_1,k_2,\cdots,k_{\eta-1},100\}=\{k_\eta\}$ such that $1\rightarrow k_\eta$. In the instance that the second case be correct, note that $\{1,k_1,k_2,\cdots,k_\eta,100\}=C\cc$ is not a closed set, implying $\exists\,j\in C$ such that $1\rightarrow j$. Therefore, in any case, $1\rightarrow j$ for some $j\in C$.\\
		
		Thus, we have a positive probability that in a finite number of moves, the process will reach the closed set $C$ and hence, the game will then become impossible to complete. That is, $N(X)$ Board is occasionally winnable.\\
		
		\noindent\textit{Necessary Part:}\\
		Let $N(X)$ Board be occasionally winnable. That is, there is a positive probability that the game will never be completed and there is  positive probability that the game will be completed in a finite number of moves. Such a case implies the process gets stuck in a set of states that cannot access state $100$. Let us have such a set $C_0$ that cannot access state $100$. Let $C\subset\bbN_{100}$ be the set of all states that can be accessed from $C_0$. Clearly, $100\notin C$ and $C$ is a closed set since all the states that could be accessed from $C$ are already included in $C$. Note that since the board is occassionally winnable (i.e., $1\rightarrow 100$), we have $1\notin C$.\\
		
		Thus, by Result \ref{downzero}, we can write the one-step transition probability matrix as $$P_{(r)}=\begin{bmatrix}P_0&Q\\\mathbb{O}&P_C\end{bmatrix}$$ with some square matrices $P_0\,(\nim p_{1,100})$ and $P_C$, the matrix of one-step transition probabilities between the states in $C$.\\
		
		Additionally, since the process gets stuck in the set of states $C$, there must exist some state $\ell\,(\in C)$ such that $1\rightarrow\ell$. Therefore, $\forall\, C_1\supset C$ such that $1,100\in C_1\cc=\bbN_{100}\backslash C_1$, $C_1\cc$ cannot be a closed set. Hence, by Result \ref{upzero}, $$P_{(r)}\neq\begin{bmatrix}P_0^*&\mathbb{O}\\Q^*&P_{C_1}\end{bmatrix}$$ with square matrices $P_0^*\,(\nim p_{1,100})$ and $P_{C_1}\,(\rhd P_C)$.
	\end{proof}

	\begin{result}[The Necessary and Sufficient Condition for Ultimately Winnable Boards]\label{ulwin}
		A \emph{Moksha-Patam} Board with one-step transition probability matrix $P$, is ultimately winnable if, and only if,
		\begin{enumerate}[(i)]
			\item there does not exist any ordering of the state space $\bbN_{100}$, such that the rearranged one-step transition probability matrix can be written as $$P_{(r)}=\begin{bmatrix}P_1&\mathbb{O}\\Q&P_{100}\end{bmatrix}$$ with square matrices $P_1\,(\nim p_{1,1})$ and $P_{100}\, (\nim p_{100,100})$, and
			\item for any rearrangement of the one-step transition probability matrix with square matrices $P_0\,(\nim p_{1,100})$ and $P'$, such that $$P_{(r)}=\begin{bmatrix}P_0&Q\\\mathbb{O}&P'\end{bmatrix},$$
			there exist square matrices $P_0^*\,(\nim p_{1,100})$ and $P''\,(\rhd P')$ such that $$P_{(r)}=\begin{bmatrix}P_0^*&\mathbb{O}\\Q'&P''\end{bmatrix}.$$
		\end{enumerate}
	\end{result}
	
	\begin{proof}
		\textit{Sufficient Part:}\\
		Let us have a \emph{Moksha-Patam} board with the two conditions listed in Result \ref{ulwin}. The first condition ensures that the board is not unwinnable [see Result \ref{unwin}]. The second condition states that if there exists some closed set $C$ not containing states $1$ and $100$ [see Result \ref{downzero}], then there must be some superset of $C$, say $C_1$, also not containing states $1$ and $100$, such that $\bbN_{100}\backslash C_1$ is a closed set [see Result \ref{upzero}], which implies $\forall j\in C_1,\: 1\centernot\rightarrow j$. That is, if there exists some closed set $C$, it is not accessible from state $1$, thus ensuring that the board is not occasionally winnable. Hence, the board must be ultimately winnable.\\
		
		\noindent\textit{Necessary Part:}\\
		Let us have an ultimately winnable \emph{Moksha-Patam} board with one-step transition probability matrix $P$. Then, $1\rightarrow 100$. Hence, $\forall A\subset\bbN_{100}$ such that $1\in A$ and $100\notin A$, $A$ cannot be a closed set. Hence, by Result \ref{upzero}, we cannot rearrange $P$ as $$P_{(r)}=\begin{bmatrix}P_1&\mathbb{O}\\Q&P_{100}\end{bmatrix}$$ with square matrices $P_1\,(\nim p_{1,1})$ and $P_{100}\,(\nim p_{100,100})$.\\
		
		Additionally, if we do have $$P_{(r)}=\begin{bmatrix}P_0&Q\\\mathbb{O}&P'\end{bmatrix},$$ with square matrices $P_0\,(\nim p_{1,100})$ and $P'$, implying we have a closed set $C$ not containing the states $1$ and $100$ [see Result \ref{downzero}], we must have $1\centernot\rightarrow j$, $\forall j\in C$. Hence, if there exist some states $i\notin C$ and $j\in C$ such that $i\rightarrow j$, then $1\centernot\rightarrow i$. Thus, if we collect all the states that can access some state in $C$, the remaining states, including states $1$ and $100$, must comprise a closed set. Hence, by Result \ref{downzero} and Result \ref{upzero}, we have the second condition. 
	\end{proof}
	
	\begin{note}
		An alternative approach to the proof of Result \ref{ulwin} can be as follows: The first condition of Result \ref{ulwin} is the negation of the condition in Result \ref{unwin}, while the second condition of the same is a negation of the conditions in Result \ref{occwin}. Hence, the conditions in this result ensure that the board is neither unwinnable nor occasionally winnable (that is, the board is ultimately winnable), and vice versa.
	\end{note}

	\section{Identifying Winnability}
	
	While we have established the mathematical criteria for the three classifications, identifying the classification of a particular board can often be a matter of careful observation, rather than working out the one-step probability matrix associated with the board. For example, we can observe the particularly curious presence of the chute-barrier in the construction of the $6(U)$ Board and the $14(G_0)$ Board, both of which are established to be unwinnable. Additionally, the occasionally winnable $7(\Delta)$ Board also has a chute-barrier, as does the $8(\Xi)$ Board. And while it may seem obvious that the presence of a chute-barrier will hinder the progression of a game towards completion, when is this hindrance sufficient to make the board unwinnable? Is it sufficient to make the board occasionally winnable? Is it necessary? We answer these questions in this section.
	
	\begin{result}[Presence of Chute-Barrier - I]\label{chbarr}
		A necessary criteria for a board to be unwinnable is the presence of a chute-barrier. In other words, it is impossible to design an unwinnable board without the use of a chute-barrier.
	\end{result}
	
	\begin{proof}
		If possible, let there exists an unwinnable board without a chute-barrier, i.e., at most five consecutive cells are the starting cells of chutes.\\
		
		Let $k=\max(\{\lambda\in\bbN_{100}:1\rightarrow\lambda\})$, i.e., $k$ is the highest cell accessible from $1$. Since the board is unwinnable, $k\neq100$.\\
		\newline
		\textsf{Case I:} $94\leq k\leq99$.\\
		In this case, state $k$ can access state $100$ in one step with probability $\frac16$, pertaining to the die outcome associated with the value $100-k$. Hence, $k\rightarrow100$. Since $1\rightarrow k$ and $k\rightarrow100$, we have $1\rightarrow100$, which implies that $k=100$.\\
		\newline
		\textsf{Case II:} $2\leq k\leq93$.\\
		In this case, at least one of the cells $k+1,k+2,\ldots,k+6$ is not the starting of a chute. Let this cell be $\lambda'$. Thus, it can either be a \emph{non-component} cell or the starting of a ladder, which ends at cell $\ell>\lambda'$. Thus, either $k\rightarrow\lambda'$ or $k\rightarrow\ell$, which implies $1\rightarrow \lambda'>k$ or $1\rightarrow\ell>k$. Either way, $$k\neq\max(\{\lambda:1\rightarrow\lambda\})).$$
		
		\noindent Hence, in both the cases, the result is established by contradiction.
	\end{proof}
	
	By using Result \ref{chbarr}, we can say that any board which does not have at least six consecutive chutes, can be won with a positive probability. This, however, is only a necessary criteria, and does not help us identify all unwinnable boards with certainty.\\
	
	\begin{result}[Presence of Chute-Barrier - II]\label{chbar}
		A necessary criteria for a board to be occasionally winnable is the presence of a chute-barrier. In other words, it is impossible to design an occasionally winnable board without the use of a chute-barrier.
	\end{result}
	
	\begin{proof}
		We have already established the role of a closed set in an occasionally winnable board in the proof of Result \ref{occwin}. Hence, it is sufficient to show that it is impossible to have a closed set (not containing state $100$ and accessible from state $1$), without the presence of a chute-barrier on the board.\\
		
		Let us assume that it is possible to have a closed set which is not containing state $100$ and is accessible from state $1$, without the presence of a chute-barrier on the board, and that we have one such closed set, $C$. Let $k$ be the largest state in $C$ that is accessible from state $1$. Now, since there are no chute-barriers, at least one of the states $k+1,k+2,\ldots,k+6$ must be non-component or have the entrance of a ladder. Hence, $\exists\, \ell>k$ which is accessible from state $k$, and hence, also accessible from state $1$. This state $\ell$ cannot belong to $C$, since $k$ is the largest state in $C$ that is accessible from state $1$. However, since $k\rightarrow l$, $C$ cannot be closed.\\
		
		Thus, the result is established by contradiction. 
	\end{proof}
	
	Note that, this line of reasoning could also be used to directly establish that positive probability of not completing the game is only possible in the presence of closed sets and closed sets cannot exist without the presence of chute-barriers, thus establishing both Result \ref{chbarr} and Result \ref{chbar} together.\\
	
	We can see from Result \ref{chbarr} and Result \ref{chbar}, that the absence of a chute-barrier implies the board is neither unwinnable nor occasionally winnable. Thus, we have our obvious criteria for an ultimately winnable board as stated below.
	
	\begin{result}[Absence of Chute-Barrier]\label{suffch}
		The lack of a chute-barrier is a sufficient criteria for a board to be ultimately winnable. In other words, if a board does not have a chute-barrier, it can almost surely be completed in a finite number of moves.
	\end{result}

	Thus, by Result \ref{suffch}, we can see that the $10(\alpha)$ Board is ultimately winnable.\\
	
	Apart from the presence or absence of chute-barriers, there are other tell-tales that may help us identify the presence of closed sets, and since it is the presence and accessibility of closed sets that decides the classification of a particular board, it may be useful to study such tell-tales. In order to do so, we make use of the following result, a detailed discussion upon which can be found in Section 3.9 of Medhi[1982]\cite{medhi}.
	
	\begin{result}\label{SPDclosed}
		$C$ is a closed class of a discrete time-homogeneous Markov Chain with state space $S$, or equivalently $\bbN_{|S|}=\{1,2,\ldots,|S|\}$, if, and only if, there exists a stationary probability distribution $\bm{\pi}=\begin{pmatrix}\pi_1&\pi_2&\cdots&\pi_{|S|}\end{pmatrix}$ such that $\forall\, i\notin C$, $\pi_i=0$ and $\forall\, i\in C$, $\pi_i\neq0$.
	\end{result}

	We can now evaluate the stationary probability distributions of the Markov Chain of a \emph{Moksha-Patam} board to detect the presence of closed sets. Of course, to evaluate the stationary probability distribution of a Markov Chain with a hundred states is a task easier said than done. However, programming languages which are built to handle large matrix operations can be of help in this regard. In order to find the stationary probability distributions, we simply find the eigenvalues of the transpose of the one-step transition probability matrix, and obtain the eigenvectors associated with the eigenvalue $1$. The transpose of these eigenvectors are the required stationary distributions.\\
	
	One stationary distribution we will always have, for any given \emph{Moksha-Patam} board, is $$\Pi=\begin{pmatrix}0&0&\cdots&0&1\end{pmatrix},$$ since for any \emph{Moksha-Patam} board, the last row of the one-step transition probability matrix is given by $$
	\begin{blockarray}{cccccccccc}
		&\textcolor{dgrey}{1}&\textcolor{dgrey}{\cdots}&\textcolor{dgrey}{94}&\textcolor{dgrey}{95}&\textcolor{dgrey}{96}&\textcolor{dgrey}{97}&\textcolor{dgrey}{98}&\textcolor{dgrey}{99}&\textcolor{dgrey}{100}\\
		\begin{block}{c[ccccccccc]}
			\textcolor{dgrey}{100}&0&\cdots&0&0&0&0&0&0&1\\
		\end{block}
	\end{blockarray}.$$
	
	And if we do find some other stationary probability distribution, we can be certain that there exists some closed set other than $\{100\}$. Of course, this does not tell us much. However, the absence of any other stationary probability distribution besides $\Pi$, gives us a certainty that there is no closed set other than $\{100\}$, and hence, the board is ultimately winnable.\\
	
	Below we give one stronger sufficient condition for ultimately winnable board.
	
	\begin{result}
		A \emph{Moksha-Patam} board, the Markov Chain associated with which has only one stationary distribution $$\Pi=\begin{pmatrix}0&0&\cdots&0&1\end{pmatrix}$$ is an ultimately winnable board.
	\end{result}
	To further narrow down these possibilities, we must first establish a few results.
	
	\begin{result}\label{escape}
		Consider a \emph{Moksha-Patam} board with a chute-barrier containing $\lambda$ chutes, starting from state $M+1$ ($=m_1$, say). Let the chute with the lowest exit in the chute-barrier end at state $m_2$, and let the chute with the lowest exit amongst the ones starting at some state belonging to $\{m_2+1,m_2+2,\dots,m_1-1\}$, end at state $m_3$, and further let the chute with the lowest exit amongst the ones starting at some state belonging to $\{m_3+1,m_3+2,\dots,m_2-1\}$, end at state $m_4$, and so on. Let the minimum of the set $\{m_i:i\in\bbN\}$ be $m$. Then, $C=\{m,m+1,\ldots,M\}$ is not closed if, and only if, there exists a ladder which starts on state $\ell\:(\in C)$ and ends on state $L\,(>M+\lambda)$.
	\end{result}
	\begin{proof}
		Observe that the sequence $\{m_i:i\in\bbN\}$ is eventually constant, since it is decreasing, bounded below by $1$ and defined on a discrete domain. Specifically, whenever $m_i=m_{i+1}$, we will have the value $m=m_i$. We prove the necessary part only, since that of the sufficient part is simple.\\
		\\
		\textit{Necessary Part:}\\
		We shall prove this part by its contrapositive statement. Let there exists no ladder which starts on state $\ell\,(\in C)$ and ends on state $L\,(>M+\lambda)$.\\
		
		Now, the only way to access a relatively lower state from the current position is via a chute. But it is given that the lowest state that can be reached from $C$ via chutes is $m$. And since there are no ladders from $C$ to some state not belonging to $C$, the largest state that can be reached from within $C$ is $M$ [see proof of Result \ref{chbarr} for a detailed discussion of the impact of chute-barrier]. Hence, no state outside $C$ can be accessed from any state in $C$. Thus, $C$ is closed.
	\end{proof}
	
	\begin{defn}[Escape Ladder]
		The ladder with the property described in Result \ref{escape} shall be known as the \emph{escape ladder} of closed set $C$.
	\end{defn}

	\begin{defn}[Ladder-Pass]
		Five or less consecutive ladders from states lesser than the smallest state of a closed set $C$ to states higher than the largest state in $C$, shall be known as a \emph{ladder-pass} (abbreviated as LP, as and when required) of $C$.
	\end{defn}

	\begin{defn}[Ladder-Bridge]
		Six or more consecutive ladders from states lesser than the smallest state of a closed set $C$ to states higher than the largest state in $C$, shall be known as a \emph{ladder-bridge} of $C$.
	\end{defn}

	\begin{defn}[Trapper]
		Let $C$ be a closed set. Then a component whose entrance does not belong to $C$ is called a \emph{trapper} of $C$ if its exit can access a cell inside $C$ with probability $1$.
	\end{defn}

	In the same manner of reasoning as we applied for chute-barriers, it can be seen that a ladder-bridge helps `jump' across a closed set no matter what the player's die throw indicates. Generally, in such scenarios, the existence of the closed set does not alter the winnability of the board, for the closed set is not accessible from state $1$. Such boards are ultimately winnable.\\
	
	Likewise, a ladder-pass allows some positive probability to the player to `jump' across the closed set, however, it is not guaranteed since the closed set is still accessible from state $1$. In such scenarios, we may have an occasionally winnable board.\\
	
	It is also possible that there exists a trapper with its entrance either less than the smallest entrance of the ladder-bridge or more than the smallest exit of the ladder-bridge. This trapper essentially leads into the closed set, hence making the combined collection of components act in the manner of a ladder-pass (providing positive probability for both never reaching the closed set and getting trapped in it). This trapper of a closed set can very well be a ladder-pass of some lower closed set. If there be six or more consecutive trappers which start from states either below the ladder-bridge or above the maximum exit of the ladder-bridge, then they cancel out the impact of the existence of a ladder-bridge entirely. Hence rises the idea of a \emph{functional} ladder-bridge (abbreviated as FLB, as and when required), where the combined effect of all these components does not interfere with the role of a ladder-bridge to bypass the closed set with probability $1$.
	
	\begin{defn}[Functional Ladder Bridge]
		The ladder-bridge of a closed set $C$ is said to be a \emph{functional}, if no trapper of $C$ has an entrance either lower than the smallest entrance or higher than the smallest exit of the said ladder-bridge.
	\end{defn}
		
	The presence of multiple chute-barriers without their escape ladders give rise to multiple closed classes. If we have ladder-bridges or ladder-passes across the combined closed set formed by those closed classes, it is easy to comment about the winnability of the board. However, if they all have their own ladder bridges and/or ladder passes, it becomes increasingly difficult to comment about the same. In such scenarios, if two chute-barriers do not have any component with its exit between them, then it is safe to say that the states between them cannot be accessed from state $1$, and hence, the two chute-barriers can be assumed to be one single chute-barrier from end to end.\\
		
\begin{figure}[H]
		\centering\scalebox{0.8}{
			\begin{tikzpicture}[node distance=2cm]
				\node (start) [startstop] {Start};
				\node (dec1) [decision, below of=start, yshift=-1cm] {Is\\ \emph{Chute-Barrier}\\ present?};
				\node (dec2) [decision, below of=dec1, yshift=-2.3cm] {Is SPD\\ other than $\Pi$\\ present?};
				\draw [arrow] (start) -- (dec1);
				\draw [arrow] (dec1) -- node[anchor=east] {Yes} (dec2);
				\node (out1) [io, right of=dec1, xshift=4cm] {Ultimately Winnable Board};
				\node (stop) [startstop, right of=out1, xshift=2cm, yshift=-7cm] {Stop};
				\draw [arrow] (dec1) -- node[anchor=south] {No} (out1);
				\draw [arrow] (out1) -| (stop);
				\draw [arrow] (dec2.east) -| node[anchor=south east] {No} (out1.south west);
				\node (pro) [process, below of=dec2, yshift=-1cm] {Identify all closed sets other than $\{100\}$};
				\draw [arrow] (dec2.south) -- node[anchor=east] {Yes} (pro);
				\node (dec3) [decision2, below of=pro, yshift=-0.3cm] {Does every closed\\ set have FLB?};
				\node (dec34) [decision2, below of=dec3, yshift=-1cm] {Do all the closed sets\\ without FLB have LP?};
				\draw [arrow] (pro) -- (dec3);
				\draw [arrow] (dec3.east) -| node[anchor=south east] {Yes} (out1);
				\node (out2) [io, right of=dec34, xshift=4cm] {Occasionally Winnable Board};
				\node (out3) [io, below of=out2] {Unwinnable Board};
				\draw [arrow] (dec34) -- node[anchor=south] {Yes} (out2);
				\draw [arrow] (dec3) -- node[anchor=east] {No} (dec34);
				\draw [arrow] (dec34) |- node[anchor=south east] {No} (out3);
				\draw [arrow] (out3) -| (stop.south);
				\draw [arrow] (out2) -| (stop.south);
		\end{tikzpicture}}
		\caption{An algorithm to identify the winnability of \emph{Moksha-Patam} boards.}\label{algo}
	\end{figure}
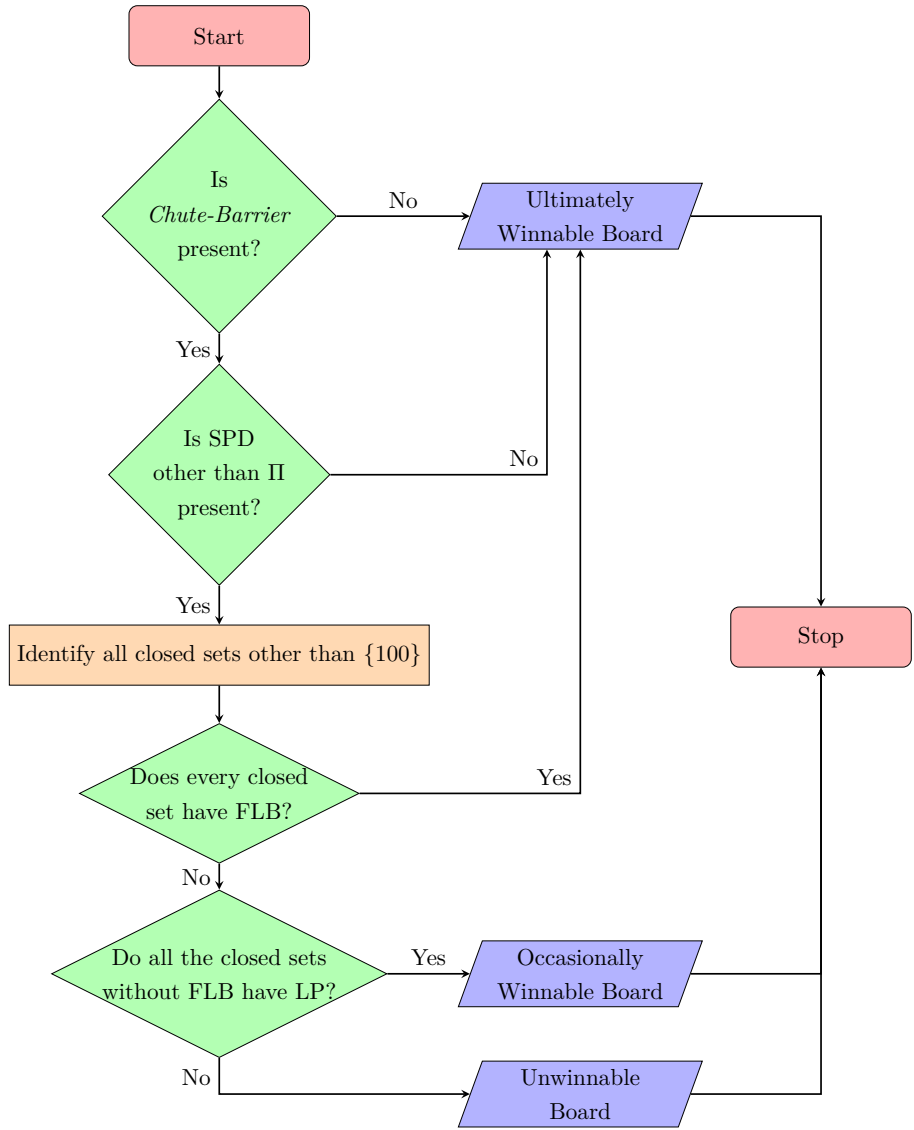

	With these observations and results in mind, we may now attempt to draw an algorithm as shown in Figure \ref{algo}.\\
	
	If stationary probability distributions (abbreviated as `SPD', as per convenience) other than $\Pi$ exist, then identification of the closed sets can be easily done by dint of Result \ref{SPDclosed}. Here, it is worth noting that the second decision criteria in the algorithm in Figure \ref{algo} asking ``Is SPD other than $\Pi$ present?'' can be replaced by the decision criteria, ``Is there a chute-barrier without an escape ladder?''. This is due to Result \ref{SPDclosed} and Result \ref{escape} combined. In that case, the identification of closed sets must be done by making use of Result \ref{escape}.\\
	
	By this algorithm, we can comment that the $8(\Xi)$ Board is ultimately winnable.
	
	\section{Too Many Boards!}
	With all the classification we have done and explored in the previous sections, it might be valid to wonder how many different \emph{Moksha-Patam} boards are possible. Interestingly, that number easily exceeds the number of atoms in the observable universe.\cite{baker}.\\
	
	\subsection{Shared Exits}
	It is difficult to find the exact number of boards possible if components are allowed to share exits. We shall obtain some crude bounds for such a count, say $\mathbb{T}$, before we obtain the exact number of possible boards with no component having shared exits. In order to do so, consider a board with $N$ components and $n(\leq N)$ distinct exits. Thus, we can construct such a board by first choosing $N+n$ cells ($N$ entrances and $n$ exits) from $98$ cells (excluding $1$ and $100$) in $\binom{98}{N+n}$ different ways. We may now choose $N$ cells from among the chosen $N+n$ cells to be the distinct entrances in $\binom{N+n}{N}$ ways. Of these $N$ components with now fixed distinct entrances, we can choose and permute $n$ components to have the $n$ distinct exits in $^NP_n=\frac{N!}{(N-n)!}$ ways. The rest of the $N-n$ components may have any of the $n$ chosen exits repeated in $n^{N-n}$ ways. This of course gives us the count $$\binom{98}{N+n}\binom{N+n}{N}\frac{N!}{(N-n)!}n^{N-n},$$
	which is a gross overestimate since the primary choice of components with distinct exits can be made in $k$ different ways if the final construction must have an exit shared by $k$ components. Thus, this count includes each board $\prod_{i=1}^n N_i$ times where $N_i$ is the number of components sharing the $i^\text{th}$ exit, with $\sum_{i=1}^n N_i=N$.\\
	
	To get a better estimate, we must divide separate terms in the count by the appropriate values of $\prod_{i=1}^n N_i$ as determined by the structure of each individual board. But using the AM-GM inequality, we have $\prod_{i=1}^n N_i\leq \left(\frac{N}{n}\right)^n$, which gives us the following crude underestimate of the required count: $$\frac{\binom{98}{N+n}\binom{N+n}{N}\frac{N!}{(N-n)!}n^{N-n}}{\left(\frac{N}{n}\right)^n}$$
	
	Intentionally ignoring the one possible board with $N=0$ for easier calculation, and summing these estimates over $n=1,2,\ldots,N$ and $N=1,2,\ldots,49$, we have the crude bounds of the total number of \emph{Moksha-Patam} boards, $\mathbb{T}$ as $$1.3135349305\times10^{127}\leq\mathbb{T}\leq1.2801985919\times10^{134}$$
	
	\subsection{No Shared Exits}	
	A fairly exact count can be obtained if we consider only the boards which do not have shared exits. To obtain this count, say $T$, we consider the total number of boards with $N$ components and $N$ distinct exits. Hence, the total number of such boards with $N$ components is given by $$\binom{98}{2N}\binom{2N}{N}N!.$$
	
	\noindent For $N=0$, of course, we have only one board, the $0$ Board.\\
	For $N=1$, the number quickly rises to $9506$. That is, there are $9506$ possible boards with only one component.\\
	For $N=2$, this number dramatically increases to $4,33,47,360$. With only two components the total number of possible boards is more than four crores!\\
	
	Hence, the total number of all possible \emph{Moksha-Patam} boards (without shared exits) is
	\begin{equation*}
		\begin{split}
			T=\sum_{N=0}^{49}\binom{98}{2N}\binom{2N}{N}N!=\sum_{N=0}^{49}\frac{98!}{(98-2N)!N!}
			&\approx7.6432896116\times10^{93}.
		\end{split}
	\end{equation*}
	$$$$
	
	Out of these total number of boards if we want to estimate the proportion of boards which are ultimately winnable (since those are the boards one would normally want to play on), then we can focus on the number of boards with chute-barriers [see Result \ref{suffch}].\\
	
	To have a chute-barrier (with $6$ chutes) in a board, we can choose the entrance of its lowest chute from amongst $\{8,9,\ldots,94\}$, since we need at least six cells below it to accommodate the exits and at least five cells above it to accommodate the consecutive entrances. Let this chosen cell be $M$. This fixes five other chutes with entrances at $M+1$, $M+2$, $M+3$, $M+4$ and $M+5$. Note that all cases of more than six chutes in the chute-barrier will be counted, in the next paragraph, when we count the possible ways to place the other components. Hence, we must now count the number ways we can choose the exits of the chutes with entrances $M, M+1, \ldots, M+5$. We can choose six cells for this purpose from $\{2,3,\ldots,M-1\}$ in $\binom{M-2}{6}$ ways and then arrange them in $6!$ ways for demarcating which exit corresponds to which chute. Hence, one such chute-barrier (with $6$ chutes) can be placed in the following number of ways: $$\sum_{M=8}^{94}\binom{M-2}{6}\times6!\approx6.8210081597\times10^{12}.$$
	
	This count can now be used as we generate boards with $N$ components which have at least one chute-barrier. We can choose $2N-12$ cells from the remaining $86$ cells in $\binom{86}{2N-12}$ ways, as entrances and exits for the remaining $N-6$ components. Out of these $2N-12$ cells, we can choose $N-6$ cells as entrances in $\binom{2N-12}{N-6}$ ways. The remaining $N-6$ cells can be arranged as exits in $(N-6)!$ ways. Hence, we can construct such an $N(X)$ board in $$\left[\sum_{M=8}^{94}\binom{M-2}{6}\times6!\right]\binom{86}{2N-12}\binom{2N-12}{N-6}(N-6)!$$ ways. However, if a board has $\ell$ chute-barriers, we could have generated the same board by starting with any one of the $\ell$ chute-barriers. Hence, this count overestimates the number of boards with more than one chute-barriers.\\
	
	Therefore, if we define $C$ as the total number of boards with chute-barriers, then $$C<\sum_{N=6}^{49}\left[\sum_{M=8}^{94}\binom{M-2}{6}\times6!\right]\binom{86}{2N-12}\binom{2N-12}{N-6}(N-6)!\approx8.6769698734\times10^{92}.$$ This gives us the following lower-bound of the number of boards without a chute-barrier: $$T-C>6.7755926243\times10^{93}.$$
	
	Since, not having a chute-barrier is a sufficient condition for an ultimately winnable board [see Result \ref{suffch}], we can say that there are at least $6.7755926243\times10^{93}$ ultimately winnable boards. Hence, if we randomly choose a \emph{Moksha-Patam} board, there is more than $88.6\%$ chance that the board will be ultimately winnable.\\
	
	In a more realistic setting where we consider a maximum of $20$ components to avoid unnecessary clutter on the board, the total number of such boards is $$\sum_{N=0}^{20}\frac{98!}{(98-2N)!N!}\approx1.7\times10^{57},$$ and the total number of such boards with chute-barriers is less than $$\sum_{N=6}^{20}\left[\sum_{M=8}^{94}\binom{M-2}{6}\times6!\right]\binom{86}{2N-12}\binom{2N-12}{N-6}(N-6)!\approx8.1\times10^{53}.$$ Hence, here we observe that there is more than $99.95\%$ chance that the board will be ultimately winnable.\\
	
	Since only ultimately winnable boards represent absorbing Markov Chains under the assumption of no shared exits, this result probably justifies the generally false belief of many, including Tun\cite{tun} and Cheteyan et. al.\cite{cheteyan} that the game of \emph{Moksha-Patam} represents an absorbing Markov Chain with the final square being the absorbing state. Of course, unwinnable boards and occasionally winnable boards can also represent absorbing Markov Chains, provided the closed classes accessible from state $1$ be singleton sets. For example, the $6\left(\begin{bmatrix}51&52&53&54&55&56\\50&50&50&50&50&50\end{bmatrix}\right)$ Board represents an absorbing Markov Chain, even though it is unwinnable.\\
	

	\bibliographystyle{plain}
	\bibliography{biblio}
\end{document}